\input{style/arxiv-general.cfg}
\documentclass[MSNbibl,number,citesort,dvips]{arxbj}
\makeatletter
   \@ifpackageloaded{graphicx}{}{\usepackage{graphicx}}
\makeatother
\usepackage{upgreek}
\usepackage{graphicx}

\aid{0}
\volume{21}
\issue{2}
\pubyear{2015}
\firstpage{1067}
\lastpage{1088}
\doi{10.3150/14-BEJ598} 

\makeatletter
\newcommand{\rrvert}{\vert}
\newcommand{\llvert}{\vert}
\newcommand{\fraca}[2]{{#1}/{#2}}
\newcommand{\eqref}[1]{(\ref{#1})}
\newremark{remark}{Remark}
\newtheorem{theorem}{Theorem}
\newtheorem{lemma}{Lemma}
\newtheorem{proposition}{Proposition}
\newproclaim{assumption}{Standing Assumption}
\newtheorem{coro}{Corollary}

\newcommand{\R}{\mathbb{R}}
\newcommand{\Rp}{\mathbb{R}_{+}}
\newcommand{\N}{\mathbb{N}}
\newcommand{\PP}{\mathbb{P}}
\newcommand{\EE}{\mathbb{E}}
\newcommand{\dd}{\mathrm{d}}
\newcommand{\FF}{\mathcal{F}}
\newcommand{\Wtil}{\tilde{W}}
\newcommand{\Ftil}{\tilde{\mathcal{F}}}
\newcommand{\supp}{\operatorname{supp}}
\newcommand{\eps}{\epsilon}
\newcommand{\sgn}{\operatorname{sgn}}
\renewcommand{\epsilon}{\varepsilon}
\renewcommand{\pi}{\uppi}
\newcommand{\cF}{\mathcal{F}}
\newcommand{\cL}{\mathcal{L}}
\newcommand{\cR}{\mathcal{R}}
\newcommand{\dto}{\downarrow}
\renewcommand{\phi}{\varphi}
\makeatother

\begin{document}
\begin{frontmatter}

\title{Finite, integrable and bounded time embeddings for diffusions}
\runtitle{Finite, integrable and bounded time embeddings for diffusions}

\begin{aug}
\author[A]{\inits{S.}\fnms{Stefan}~\snm{Ankirchner}\corref{}\thanksref{A,e1}\ead[label=e1,mark]{ankirchner@hcm.uni-bonn.de}},
\author[B]{\inits{D.}\fnms{David}~\snm{Hobson}\thanksref{B}\ead[label=e2]{D.Hobson@warwick.ac.uk}} \and
\author[A]{\inits{P.}\fnms{Philipp}~\snm{Strack}\thanksref{A,e3}\ead[label=e3,mark]{pstrack@uni-bonn.de}}
\address[A]{Hausdorff Center for Mathematics and Institute for Applied Mathematics,
University of Bonn, Endenicher Allee 60, 53115 Bonn, Germany.\\ \printead{e1,e3}}
\address[B]{Department of Statistics, University of Warwick, Coventry,
CV4 7AL, UK.\\ \printead{e2}}
\end{aug}

\received{\smonth{8} \syear{2013}}
\revised{\smonth{1} \syear{2014}}

%
\begin{abstract}
We solve the Skorokhod embedding problem (SEP) for a general \textit{time-homogeneous diffusion} $X$: given a distribution $\rho$, we
construct a
stopping time $\tau$ such that the stopped process $X_\tau$ has the
distribution $\rho$. Our solution method makes use of martingale
representations (in a similar way to Bass
(In \textit{Seminar on Probability {XVII}}. \textit{Lecture Notes in Math.}
\textbf{784}
(1983) 221--224 Springer)
who solves the SEP
for Brownian motion) and draws on law uniqueness of weak solutions of SDEs.

Then we ask if there exist solutions of the SEP which are respectively finite
almost surely, integrable or bounded, and when does our proposed construction
have these properties. We provide conditions that guarantee existence
of \textit{finite time} solutions. Then, we fully characterize the distributions
that can
be embedded with \textit{integrable} stopping times. Finally, we derive
necessary, respectively sufficient, conditions under which there exists
a \textit{bounded} embedding.
\end{abstract}

%

\end{frontmatter}

\section*{Introduction}
\label{sec:intro}

Let $X$ be a one-dimensional time-homogeneous diffusion, and let $\rho
$ be a
probability measure on $\R$. The Skorokhod embedding problem (SEP) for
$\rho$
in $X$ is to find a stopping time $\tau$ such that $X_\tau\sim\rho$. Our
main goals in this article are firstly to construct a solution of the
Skorokhod embedding problem, and secondly to discuss when does there
exist a
solution which is finite, integrable or bounded in time, and when does our
construction have these properties.

Our construction is based on the observation that we can remove the
drift of
the time-homogeneous diffusion by changing the space variable via a
\textit{scale function}. We can thus simplify the embedding problem to the case where
$X$ is a local martingale diffusion. We then consider a random variable that
has the distribution we want to embed and that can be represented as a
Brownian martingale $N$ on the time interval $[0,1]$. Further, we set
up an
ODE that uniquely determines a time-change along every path of $X$. We then
show, by drawing on a result of uniqueness in law for weak solutions of SDEs,
that the time-changed diffusion has the same distribution as the martingale
$N$. Thus the time-change provides a solution of the SEP.

Our solution is a generalization of Bass's solution of the SEP for
\textit{Brownian motion} (see \cite{Bass1983}). Bass
also starts with the martingale
representation of a random variable with the given distribution. By changing
the martingale's clock, he obtains a Brownian motion and an associated
embedding stopping time. The time-change is governed by an ODE, a
special case
of our ODE, which establishes an analytic link between Brownian paths and
the embedding stopping time. This link yields embedding stopping
times for \textit{arbitrary} Brownian motions.

Now consider properties of solutions of the SEP.
As is well known from the literature, whether a distribution is embeddable
into the diffusion $X$ depends on the relation between the support of the
distribution and the state space of $X$ and the relation between the initial
value $X_0$ and the first moment of the distribution. We include in our
analysis a general discussion of sufficient and necessary conditions
for the existence of \emph{finite} embedding stopping
times, with particular reference to our proposed construction.

Next, we
fully determine the collection
of distributions that can be embedded in $X$ with \textit{integrable} stopping
times. The associated conditions involve an integrability condition on the
target distribution which makes use of a function that also
appears in Feller's test for explosions (see, e.g., \cite{KS}).



Finally, we address the question of whether a distribution can be
embedded in
\textit{bounded} time. Recall that the Root solution (\cite{Root1969}) of the
SEP has the property that it minimises
$\EE[(\tau- t)^+]$ uniformly in $t$. The Root solution $\tau_R$ is
of the
form \( \tau_R = \inf\{ t \dvtx  (X_t, t) \in\cR_{\rho} \} \) where
$\cR_\rho$ is
a `barrier'; in particular $\cR_{\rho} = \{ (x,t) \dvtx  t \geq\beta
_\rho(x) \}
\subseteq\R\times\Rp$ for some suitably regular function $\beta
_\rho$
depending on the target law. Hence, a necessary and sufficient
condition for
there to be an embedding $\tau$ with $\tau\leq T$ is that $\beta
_\rho( \cdot
) \leq T$. However, the Root barrier is non-constructive and difficult to
analyse (though for some recent progress in this direction see Cox and
Wang \cite{CoxWang2011} and Oberhauser and Dos Reis \cite{OberhauserDosReis}).
For this reason, instead of searching for a single set of necessary
\textit{and}
sufficient conditions we limit ourselves to finding separate sets of necessary
conditions and sufficient conditions.

Our original motivation in developing a solution of the SEP for diffusions
was to study bounded stopping times with the aim of providing simple
sufficient conditions for the existence of a bounded embedding.
The boundedness (finiteness) of an embedding is an important property
of the embedding used to solve the gambling in contests problem of Seel and
Strack \cite{SeelStrack2008}, and is also relevant in the model-independent
pricing of variance swaps, see Carr and Lee \cite{CarrLee},
Hobson \cite{Hobson2011} and Cox and Wang \cite{CoxWang2011}.

Consider for a moment the case where $X$ is a real-valued Brownian
motion, null
at 0. Then it is possible to embed any target probability measure $\rho
$ in
$X$. Moreover, $\rho$ can be embedded in integrable time if and only
if $\rho$
is centred and in $L^2$, and then $\EE[\tau] = \int x^2 \rho(\mathrm{d}x)$.
The case of
embeddings in bounded time is more subtle.
Clearly a necessary condition for there to exist an embedding $\tau$
of $\rho$
in $X$ such that $\tau\leq1$ is that $\rho$ is smaller that $\mu_G$ in
convex order, where $\mu_G$ is the law of a standard Gaussian. But
this is not
sufficient. Let $\mu_{\pm a}$ be the uniform measure on $\{-a, +a \}$. Then
$\mu_{\pm a}$ is smaller than $\mu_G$ in convex order if and only if
$a \leq
\sqrt{2/\pi}$. But any embedding $\tau$ of $\mu_{\pm a}$ has $\tau
\geq\min
\{ u \dvtx  |B_{u}| \geq a\}$, and thus does not satisfy $\tau\leq T$ for
any $T$.
Hence, we would like to find sufficient conditions on $\rho$ such that there
exists $\tau\leq T$ with $X_\tau\sim\rho$. The case where $X$ is Brownian
motion, possibly with drift, was considered in Ankirchner and
Strack \cite{as11}. Here we consider general time-homogeneous diffusions.

The paper is organized as follows. In Section~\ref{sec:mg}, we
describe our
solution method of the SEP for a diffusion \textit{without} drift. In this
section, we assume that the initial value of the diffusion coincides
with the
first moment of the distribution to embed (the centred case). In the following
Section~\ref{non-c case}, we briefly explain how to construct
solutions if the
first moment does \textit{not} match the initial value (the
non-centred case). In
Section~\ref{sec finite e}, we collect some general conditions which guarantee
that a distribution can be embedded into $X$ in finite time. We then consider
integrable embeddings in Section~\ref{sec int emb}. Distinguishing between
the centred and non-centred case, we provide sufficient and necessary
conditions for the existence of integrable solutions of the SEP.
Section~\ref{sec:bte} discusses bounded embeddings. In Section~\ref{sec:diff}, we
explain how one can reduce the SEP for diffusions with drift to the case
without drift. Finally, in Section~\ref{sec:eg} we illustrate our
results with
some examples.

\section{The martingale case}
\label{sec:mg}

We will argue in Section~\ref{sec:diff} below that the problem of
interest can
be reduced to the case in which the process is a continuous local martingale.
In this section, we describe a generalisation of the Bass \cite{Bass1983}
solution of the SEP. The Bass solution is an embedding of $\nu$ in Brownian
motion: we consider embeddings in a local martingale diffusion which
may be
thought of as time-changed Brownian motion.

Consider the time-homogeneous local martingale diffusion $M$, where $M$ solves
%
\begin{equation}
\label{sde} \mathrm{d}M _s = \eta(M_s)\,\mathrm{d}W
_s, \qquad\mbox{with } M_0=m;
\end{equation}
here
$m \in\R$ and $\eta\dvtx \R\to\R_{+}$ is Borel-measurable. We assume
that the
set of points $x\in\R$ with $\eta(x) = 0$ coincides with the set of points
where $\frac{1}{\eta}$ is \textit{not} locally square integrable.
Then a result
by Engelbert and Schmidt implies that the SDE \eqref{sde} possesses a weak
solution that is unique in law (see, e.g., Theorem~5.5.4 in \cite
{KS}). We define
$l = \sup\{x \leq m \mid \eta(x) = 0\}$ and $r = \inf\{x \geq m \mid \eta
(x) =
0\}$ so that $-\infty\leq l \leq r \leq\infty$ (to exclude
trivialities we
assume $l < m < r$) and for $x \in\R$,
%
\begin{equation}
\label{eqn:qdef} q(x) = \int_m^x \,\mathrm{d}y
\int_m^y \frac{2}{\eta^2(z)} \,\mathrm{d}z.
\end{equation}
By our assumption on $\eta$, $q$ is infinite on $(-\infty, l)$ and
$(r, \infty)$.

\begin{remark}
\label{rem:felkot}
By Feller's test, $\PP[\inf_{s \le t} M_s
\le l] = 0$ for one, and then every, $t>0$ if and only if $q(l+) = \lim_{x
\downarrow l} q(x) = \infty$. Similarly, $\PP[\sup_{s \le t} M_s \ge
r] = 0$
if and only if $q(r-) = \infty$ (see, e.g., Theorem~5.5.29 in \cite{KS}).
Further, by results of Kotani \cite{Kotani2006}, the local martingale
$M$ is a
martingale \emph{provided} either $-\infty< l$ or $\int_{-\infty} |x|
\eta(x)^{-2} \,\mathrm{d}x  = \infty$ \emph{and} either $r<\infty$ or $\int^\infty x
\eta(x)^{-2} \,\mathrm{d}x  = \infty$.
\end{remark}

Note that our assumption that $\frac{1}{\eta}$ is not locally square
integrable at $l$ and $r$ implies that $l$ and $r$ are absorbing
boundaries if
they can be reached in finite time. Then without loss of generality we may
assume that $\eta= 0$ on $(-\infty,l)$ and $(r, \infty)$ and $\eta$ is
positive on $(l,r)$.\looseness=-1

We want to embed a non-Dirac probability measure $\nu$ with $ \int x
\,\mathrm{d}\nu(x) =
m$.
Let $\underline{v} = \inf\{ \supp( \nu) \}$
and $\overline{v} = \sup\{ \supp( \nu) \}$ be the extremes of the
support of
$\nu$,
and
let $F$ be the distribution function associated to the target law $\nu$.
Moreover, let $\Phi\dvtx  \R\rightarrow[ 0, 1 ]$ be the cumulative
distribution
function of the normal distribution and $\phi= \Phi'$ its density.
Define the
function $h = F^{- 1} \circ\Phi$. Let $(\Wtil_t)_{t\geq0}$ be a Brownian
motion on a filtration $\tilde{\mathbb F} = (\Ftil_t)_{t \geq0}$.
Notice that
$h(\Wtil_1)$ has the distribution $\nu$. In particular, $h(\Wtil_1)$ is
integrable and $\EE[h(\Wtil_1)]=m$.

We define the $\tilde{\mathbb F}$-martingale
$N_t = \EE[ h ( \Wtil_1 ) \mid \tilde{\cF}_t ]$ for $t \in
[0,1]$.
Notice that $N_0=m$, $N_1$ has distribution $\nu$ and $N_t = b(t,\Wtil_t)$,
where
\[
b(t,x) = \int_{\R} h ( y ) \phi_{1 - t} (x - y) \,
\mathrm{d}y = (\phi_{1-t} \star h) (x),
\]
and $\phi_v$ is the density of the normal distribution with variance $v$.

Since $\nu$ is not a Dirac measure we have that $h$ is increasing
somewhere, and hence, for all $t \in[0,1)$, the mapping $x \mapsto b
(t, x)$ is strictly increasing. Thus, we can define the inverse function
$B \dvtx  [ 0, 1] \times\R\rightarrow\R$ implicitly by
%
\begin{equation}
\label{defi b} b \bigl( t, B ( t, x ) \bigr) = x, \qquad\mbox{for all } t \in[ 0,
1 ), x \in\R;
\end{equation}
moreover we set $B(1,x) = h^{-1}(x)$.
The process $N$ solves the SDE
%
\begin{equation}
\label{sde2} dN_t = b_x\bigl(t,B(t,N_t)
\bigr) \,\mathrm{d}\Wtil_t,\qquad N_0 = \int x \,
\mathrm{d}\nu(x) = m.
\end{equation}
Define
%
\begin{equation}
\label{eq:defrR} \lambda(t,y)= \frac{b_x(t,y)}{\eta(b(t,y))},\qquad \Lambda(t,y)=\lambda
\bigl(t, B(t,y)\bigr) = \frac{b_x(t,B(t,y))}{\eta(y)}.
\end{equation}
The candidate embedding which we want to discuss is $\delta(1)$ where
$\delta$
solves
%
\begin{equation}
\label{ode} \delta'(t) = \Lambda(t, M_{\delta(t)})^2
= \frac
{b_x(t,B(t,M_{\delta(t)}))^2}{\eta(M_{\delta(t)})^2},\qquad \delta(0) = 0.
\end{equation}
Note that $\delta$ is increasing so that if $\delta$ is defined on
$[0,1)$ then
we can set $\delta(1) = \lim_{t \uparrow1} \delta(t)$.

\begin{theorem} \label{thm:odesol}
If the ODE \eqref{ode} has a
solution on
$[0,1)$
for almost all paths of $M$, then $\delta(1)$ embeds $F$ into $M$,
that is, the
law of $M_{\delta(1)}$ is $\nu$.
\end{theorem}

\begin{pf}
Let $Y_t = M_{\delta(t)}$ for all $t \in[0,1) $. By interchanging the
time-change and integration, see, for example, Proposition V.1.5 in
\cite{RY},
we get
\[
Y_t - m = \int_0^{\delta(t)}
\eta(M_s) \,\mathrm{d}W _s = \int_0^t
\eta(M_{\delta(s)}) \,\mathrm{d}W _{\delta(s)} = \int_0^t
\eta(Y_s) \,\mathrm{d}W _{\delta(s)}.
\]
Let $Z_t = \int_0^t \frac{1}{\sqrt{\delta'(s)}} \,\mathrm{d}W _{\delta(s)}$,
for $t
\in[0,1]$. Notice that $\langle Z,Z \rangle_t = \int_0^t
\frac{1}{\delta'(s)} \,\mathrm{d}\delta(s) = t$ (Proposition V.1.5 in
\cite{RY}) and then by L\'evy's characterization theorem, $Z$ is a
Brownian motion
on $[0,1]$.
Next, observe that
%
\begin{eqnarray}
\label{Ysol} Y_t - m &=& \int_0^t
\eta(Y_s) \sqrt{\delta'(s)} \,\mathrm{d}Z
_s = \int_0^t
\eta(Y_s) \Lambda(s, M_{\delta(s)})\,\mathrm{d}Z _s
\nonumber
\\[-8pt]
\\[-8pt]
&=& \int_0^t b_x
\bigl(s,B(s,Y_{s})\bigr) \,\mathrm{d}Z _s,
\nonumber
\end{eqnarray}
which shows that $Y$ solves the SDE \eqref{sde2} with $W$ replaced by $Z$;
in other words $(Y,Z)$ is a weak solution of \eqref{sde2}.

It follows directly from Lemma~2(a) in Bass \cite{Bass1983} that
$b_x(t,B(t,x))$ is Lipschitz continuous in $x$, uniformly in $t$, on compact
subsets of $[0,1) \times\R$. Therefore, the SDE \eqref{sde2} has at
most one
strong solution on $[0,1)$ and hence \eqref{sde2} is pathwise unique, from
which it follows (see, e.g., Section~5.3 in \cite{KS}) that solutions of
\eqref{sde2} are unique in law. Hence, for $t<1$, $Y_t=M_{\delta(t)}$
has the
same distribution as $N_t$, and in the limit $t$ tends to 1 we have
$N_1$ and
hence $Y_1$ has law $\nu$.
\end{pf}

\begin{remark}
Notice that the assumption that $\int x \nu(\mathrm{d}x) = m$ is crucial
for the conclusion of Theorem~\ref{thm:odesol}. Indeed, if $\int x \nu(\mathrm{d}x)
\neq m$, then $Y$ and $N$ solve the same SDE, but with \emph{different}
initial conditions. Hence, one cannot derive that $Y_1$ has the same
distribution as $N_1$.
\end{remark}

We next aim at showing that $\delta(1)$ is a stopping time with
respect to
${\mathbb F}^M = (\cF^M_t)_{t \geq0}$, the smallest filtration
containing the
filtration generated by the martingale $M$ and satisfying the usual
conditions. To this end we consider, as in \cite{Bass1983}, the ODE satisfied
by the inverse of $\delta(t)$. The ODE for the inverse is Lipschitz continuous
and hence guarantees that Picard iterations converge to a unique solution.

\begin{lemma}\label{lemma ode inv}
Let $M$ be a path of the solution of \eqref{sde}.
Then \eqref{ode} has a solution on $[0,1)$ if and only if there exists
$a \in
\Rp\cup\{\infty\}$ such that the ODE
%
\begin{equation}
\label{eq:Delta} \Delta'(s) = \frac{\eta(M_s)^2}{b_x(\Delta(s),
B(\Delta(s),M_s))^2}
\end{equation}
has a solution on $[0,a)$ with $\lim_{s \uparrow a} \Delta(s) = 1$.
\end{lemma}

\begin{pf}
Assume that there exists a solution of \eqref{ode} on $[0,1)$. Set $a =
\delta(1)$ and define $\Delta(s) = \delta^{-1}(s)$ for all $s \in
[0,a]$. Then a straightforward calculation shows that $\Delta$ satisfies
\eqref{eq:Delta}.\looseness=-1

The reverse direction can be shown similarly.
\end{pf}

\begin{remark}
If $\eta= 1$, then the ODE \eqref{eq:Delta} is the ODE (1)
of Bass'
paper \cite{Bass1983}.
\end{remark}


\begin{lemma}
Suppose the ODE \eqref{ode} has a solution on $[0,1)$ for almost all
paths of
$M$. Then $\delta(t)$ is an ${\mathbb F}^M$-stopping time, for all
$t\in[0,1]$.
\end{lemma}

\begin{pf}
Let $C$ be a compact subset of $[0,1) \times\Rp$. By
Lemma~2 of Bass \cite{Bass1983}, $b_x(t,x)$ and $B(t,x)$ are Lipschitz
continuous on $C$. Moreover, on $C$ the function $b_x$ is bounded away from
zero and bounded from above. This implies that $\frac
{1}{b_x(t,B(t,x))^2}$ is
Lipschitz continuous on $C$, too.

Define the mapping $\gamma\dvtx  (t,y) \mapsto\frac{\eta(M_t)^2}{b_x(y,
B(y,M_t))^2}$. Now let $D$ be a compact subset of $\Rp\times[0,1)$.
Then there exists an
$L \in\Rp$ such that for all $(t,y)$ and $(t, \tilde y) \in D$ we have
\[
\bigl\llvert \gamma(t,y) - \gamma(t,\tilde y)\bigr\rrvert \leq L
\eta(M_t)^2 \llvert y - \tilde y\rrvert ,
\]
that is, $\gamma$ is Lipschitz continuous in the second argument.

We define the Picard iterations $\Delta_0(t) = 0$ and for $n \geq0$,
\[
\Delta_{n+1}(t) = 1 \wedge\int_0^t
\frac{\eta(M_s)^2}{b_x(\Delta_n(s),
B(\Delta_n(s),M_s))^2} \,\mathrm{d}s.
\]
We have that $\Delta_n(t) = 1$ after the first time where
$\Delta_n$ attains $1$.
The assumptions on $\eta$ guarantee that $\int_0^s \eta(M_t)^2 \,\mathrm{d}t $
is finite, a.s.
for each $s<a$
(see, e.g., Section~5.5 in \cite{KS}).
By standard arguments, one can now show that $\Delta_n(t)$ converges to
a limit $\Delta(t)$ on $[0,a)$, where $\Delta(a) = 1$ for all
$t \ge a$. In particular, $\Delta(t)$ is $\cF^M_t$-measurable; moreover
$\Delta(t)$ solves the ODE \eqref{eq:Delta} on $[0,a)$.

Now let $t \in[0,1)$ and $u \in\Rp$. Observe that
\[
\bigl\{\delta(t) \le u\bigr\} = \bigl\{ \Delta(u) \ge t\bigr\}.
\]
The RHS is $\cF^M_u$-measurable, which implies that $\delta(t)$ is an
$(\cF^M_t)$-stopping time. The limit $\delta(1) = \lim_{t \uparrow1}
\delta(t)$ is also an $(\cF_t^M)$-stopping time.
\end{pf}

\begin{lemma}\label{aux lem 160612}
There exists a solution of \eqref{ode} on $[0,1)$ for almost all paths
of $M$ if
and only if
$\int_0^T \Lambda(t,N_t)^2 \,\mathrm{d}t  < \infty$, a.s. for all $T<1$. In this
case, $\delta(1)$
has the same
distribution as $\int_0^1
\Lambda(t,N_t)^2 \,\mathrm{d}t $.
\end{lemma}

\begin{pf}
For all $n \in\N$ let $\xi_n = 1 \wedge\inf\{t \ge0\mid \int_0^t
\Lambda(s,
M_{\delta(s)})^2 \,\mathrm{d}s  \ge n\}$ and $\zeta_n = 1 \wedge\inf\{t \ge0\mid
\int_0^t
\Lambda(s, N_s)^2 \,\mathrm{d}s  \ge n\}$. By appealing to uniqueness in law of solutions
of \eqref{sde2} one can show, as in the proof of Theorem~\ref{thm:odesol},
that $(M_{\delta(t)})_{0 \le t \le\xi_n}$ and $(N_t)_{0 \le t \le
\zeta_n}$
have the same distribution. Moreover, $\xi_n$ and $\zeta_n$ have the same
distribution, and therefore, $\lim_n \PP[\xi_n = 1] = 1$ if and only if
$\lim_n \PP[\zeta_n = 1] =1$.
\end{pf}

Recall (Monroe \cite{Monroe1972}) that a solution $\sigma$ of the SEP for
$\nu$ in $M$ is \emph{minimal} if whenever $\tau$ is a solution of the
SEP for
$\nu$ in $M$ such that $\tau\leq\sigma$ then $\tau= \sigma$
almost surely.
The following result shows that $\delta(1)$ is minimal, provided it
exists. In
particular, the Bass embedding \cite{Bass1983} is minimal.

\begin{proposition}
Suppose $\int_0^t \Lambda(s,N_s)^2 \,\mathrm{d}s  < \infty$ almost surely, for
every $t<1$, or
equivalently $\delta(t)<\infty$ almost surely for each $t<1$. Then
$\delta(1)$ is a
minimal embedding of $\nu$ in $M$.\looseness=-1
\end{proposition}

\begin{pf}
We have $(N_t)_{0 \leq t \leq1} = (\EE_t[h(\tilde{W}_1])_{0 \leq t
\leq1}$
is uniformly integrable (UI). Since, by construction $Y \stackrel{\cL
}{=} N$,
it follows that $(Y_t)_{0 \leq t \leq1}$ is UI. But $Y_t \equiv
M_{\delta(t)}$
and $M_s = W_{A_s}$ for some time-change $A$ and some Brownian motion
$W$ and
hence $(W_{A \circ\delta(s)})_{0 \leq s \leq1} = (W_{s \wedge(A
\circ
\delta)(1)})_{s \geq0}$ is UI. Monroe \cite[Theorem~3]{Monroe1972} proves
that in the Brownian case, if $\tau$ is an embedding of $\nu$ in a Brownian
motion $W$ and if $W_0 = \int x \nu(\mathrm{d}x)$ then $\tau$ is minimal if
and only if
$(W_{t \wedge\tau})_{t \geq0}$ is UI. Hence, $A_{\delta(1)}$ is
minimal for
$\nu$ in $W$. Since $A$ is increasing we can conclude that $\delta
(1)$ is
minimal for $\nu$ in $M$.
\end{pf}

\begin{theorem}
Suppose $\supp(\nu) \subseteq[l,r]$. Recall $M_0=m$ and suppose $\nu
\in L^1$
and $\int x \nu(\mathrm{d}x)=m$. Then $\delta(1)$ exists and is a minimal
embedding of
$\nu$ in $M$.
\end{theorem}

\begin{pf}
For $t<1$, $N_t \in(\underline{v},\overline{v})
\subseteq(\ell,r)$ and since $\frac{1}{\eta}$ is locally square integrable
$\int_0^t
\Lambda(s,\break N_s)^2 \,\mathrm{d}s  =
\int_0^t \frac{b_x(s, B(s,N_s))^2}{\eta(N_s)^2} \,\mathrm{d}s  < \infty$ almost surely.
Hence,
$\delta(t)$ exists and is finite for each $t<1$ and $M_{\delta(1)}$ has
law $\nu$.
\end{pf}

Suppose $\nu$ places mass outside $[\ell,r]$. Then it is clear that
it is not
possible to embed
$\nu$ in $M$ using any embedding. To see that this holds true for
$\delta(1)$,
suppose $\overline{v}>r$. Then for each $t<1$ we have $b(t, \cdot) \dvtx
\R\to
(\underline{v},\overline{v})$ and there exists a continuous function $y(t)$
such
that $b(t,y) > r$ for $y> y(t)$. Then,
$\int_0^T \lambda(t, \tilde{W}_t)^2 \,\mathrm{d}t  = \infty$ for all $T<1$ such that
$\sup_{0<s<T}
\tilde{W}_s - y(s)>0$. Since the set $\sup_{0<s<1} \tilde{W}_s -
y(s)>0$ has
positive probability, $\delta$ explodes before time 1 with positive
probability
also.

\begin{assumption}
Henceforth, we will assume that $\nu$ places no mass outside $[\ell,r]$.
\end{assumption}

Recall that we have assumed that we are given a diffusion with $M_0=m$, and
that the target measure $\nu$ satisfies $\nu\in L^1$ and $m = \int x
\nu(\mathrm{d}x)$. We call this the centred case. In the next section, we
consider what
happens if we relax this assumption.

In the case where $\nu\in L^1$ but $m \neq\int x
\nu(\mathrm{d}x)$, we introduce an embedding $\delta^*$ which involves running the
martingale $M$ until it first hits $\int x \nu(\mathrm{d}x)$ and then using the
stopping time $\delta(1)$ defined above, but for $M$ started at $\int x
\nu(\mathrm{d}x)$.

Then in subsequent sections we will ask, when does there exist a finite
(respectively \{integrable, bounded\}) embedding, and when does $\delta
(1)$ or
more generally $\delta^*$ have this property.

\section{The non-centred case}\label{non-c case}

In this section, we do not assume that $\nu\in L^1$ and that $m= \int x
\nu(\mathrm{d}x)$.

When at least one of $\int_{-\infty}
|x| \nu(\mathrm{d}x)$ and $\int^\infty x\nu(\mathrm{d}x)$ is finite we write
${\nu}^* = \int x \nu(\mathrm{d}x) \in\bar{\R}$.
Note that we assume that $\nu$ has support in the state space of $M$.

\begin{proposition}[(Pedersen and Peskir \cite{PedersenPeskir2001}, Cox and
Hobson \cite{CoxHobson2004})] \label{prop:mneqbarnu}
Suppose $-\infty< l < m <\allowbreak   r < \infty$. Then for there to be an
embedding of
$\nu$ in $M$ we must have that $\int x \nu(\mathrm{d}x) = m$. In this case $M$
is a
uniformly integrable martingale.

Suppose $-\infty=l < m < r < \infty$. Then there exists an embedding
of $\nu$
in $M$ if and only if ${\nu}^* \geq m$. Conversely, if $-\infty< l <
m < r
= \infty$ there exists an embedding of $\nu$ in $M$ if and only if
${\nu}^*
\leq m$.

Finally, suppose $-\infty=l < m < r = \infty$. Then we can embed any
distribution
$\nu$ in $M$.
\end{proposition}

\begin{pf}
In the bounded case, the fact that $M$ is a bounded local
martingale gives that it is a UI-martingale, and hence $\int x \nu(\mathrm{d}x) =
\EE[M_\tau] = M_0=m$.

For the second case, the upper bound on the state space means that $M$
is a
submartingale so that the condition $m \leq{\nu}^*$ is necessary. Then
provided $\nu\in L^1$ we can run $M$ until it first reaches ${\nu}^*
\in
[m,\infty)$. Note that $M$ hits $\nu^*$ in finite time by the
argument in
Karatzas, Shreve \cite{KS}, Section~5.5 C. Then we can embed $\nu$
using the
local martingale $M$ started from ${\nu}^*$ (using, for example, the time
$\delta(1)$ defined above, or the Az{\'e}ma--Yor construction as in
Pedersen and
Peskir \cite{PedersenPeskir2001}). If ${\nu}^*$ is infinite, then we
need a
different construction, see, for example, Cox and Hobson \cite{CoxHobson2004}.

For the final case, any distribution can be embedded in $M$. If $\nu
\in L^1$
then we can run $M$ until it hits ${\nu}^*$ and then consider an
embedding for
the local martingale started at the mean of the target distribution. If
$\nu
\notin L^1$,  then we can use the construction in \cite
{CoxHobson2004}, but not
the one in this paper.
\end{pf}

Let $H^M_{z}$ be the first hitting time of $z$ by
$M$, and more generally let $H^X_x$ be the first hitting time of $x$ by a
stochastic process $X$. Suppose $\mu\in L^1$ and let $\delta_{\nu
^*}(1)$ be
the
stopping
time $\delta(1)$ constructed
in the previous section to embed $\nu$ in the time-homogeneous diffusion
started at $M_0=\nu^*$. Then let $\delta^* = H^M_{\nu^*} +
\delta_{\nu^*}(1)$. By the results of the proposition, provided $\nu
\in L^1$
and both $\nu^* \leq m$ if $r<\infty$ and $\nu^* \geq m$ if
$l>-\infty$, then
$\delta^*$ is an embedding of $\nu$.

\section{Finite embeddings}
\label{sec finite e}

\subsection{The centred case}
Suppose $\nu\in L^1$ and $m = \int x \nu(\mathrm{d}x)$.

\begin{proposition}
\begin{enumerate}[(ii)]
\item[(i)]If $\ell>-\infty$, $M$ does not hit $\ell$ in finite time and
$\nu(\{l\})>0$
or if
$r<\infty$, $M$ does not hit $r$ in finite time and $\nu(\{r\})>0$, then
any embedding of $\nu$ has $\tau= \infty$ with positive probability.

\item[(ii)]Otherwise, either $\ell= -\infty$, or $M$ does not hit $\ell$
in finite
time
and $\nu(\{ \ell\}) = 0$ or $M$ can hit $\ell$ in finite time \emph{and}
either $r = \infty$, or $M$ does not hit $r$ in finite time and $\nu
(\{ r \})
= 0$ or $M$ can hit $r$ in finite time. Then if $\tau$ is an embedding of
$\nu$ we have that $\bar{\tau}= \tau\wedge H^M_\ell\wedge H^M_r$
is also an
embedding of $\nu$ and $\bar{\tau}$ is finite almost surely.\looseness=1
\end{enumerate}
\end{proposition}

\begin{pf}
(i) Suppose $\tau$ is any embedding of $\nu$ in $M$. Then
$\tau= \infty$ on the set where $M_{\tau} \in\{ \ell, r \}$.
Moreover, this set has positive probability by assumption.

(ii) If $\tau$ is an embedding of $\nu$, then $M_{t \wedge\tau}$ converges
almost surely, even on the set $\tau= \infty$. However, if $(\ell=
-\infty,r=\infty)$ then by the Rogozin trichotomy (see \cite
{rogozin}), $- \infty= \liminf M_t <
\limsup M_t = \infty$ and $(M_t)_{t \geq0}$ does not converge. Hence,
we must
have $\tau<\infty$.\looseness=1

Otherwise, one or both of $\{\ell, r\}$ is finite. Then
$M$ converges and so if $\tau= \infty$ then either
$M_{\tau}=\ell$ or $M_{\tau}=r$.

If $\ell$ or $r$ is finite but $M$ hits neither $\ell$ nor $r$ in
finite time,
then $\tau= \infty$ is excluded outside a set of measure zero by the
hypothesis that $\nu(\{\ell\})=0$ and $\nu(\{r\})=0$. Hence, $\tau<
\infty$
almost surely.

Finally, if $M$ can hit either $\ell$ or $r$ in finite time then it
will do so
and $\bar{\tau}= H^M_\ell\wedge\break H^M_r < \infty$.
\end{pf}

\begin{coro}
\label{cor:finite}
If there exists an embedding $\tau$ of $\nu$ in $M$ which is
finite almost surely then $\delta(1)$ is finite almost surely.
\end{coro}

\begin{pf}
If there is a finite embedding, then we must be in case (ii) of the
proposition. Then $\delta(1)\wedge H^M_\ell\wedge H^M_r$ is finite almost
surely. But also $\delta(1)\leq H^M_\ell\wedge H^M_r$ so that
$\delta(1)<\infty$ almost surely.
\end{pf}

\subsection{The non-centred case}

Suppose $\nu$ and $m$ are such that an embedding exists (recall
Proposition~\ref{prop:mneqbarnu}). Necessarily we must have that at
least one
of $\ell$ and $r$ is infinite.

Suppose $\nu\in L^1$ so that $\nu^*$ and $\delta^*$ are well
defined. Then
since $H^M_{\nu^*}$ is finite almost surely, we have that $\delta^*$
is finite
if and only if $\delta_{\nu^*}(1)$ is finite almost surely.

Then the result for the non-centred case is identical to both the proposition
and the corollary describing the results in the centred case, modulo the
substitution of $\delta^*$ for $\delta(1)$ in Corollary~\ref{cor:finite}.

\section{Integrable embeddings}
\label{sec int emb}

\subsection{The centred case}

Suppose $\nu\in L^1$ and $m = \int x \nu(\mathrm{d}x)$.

In this section, we provide an integrability condition on $\nu$ that guarantees
that
\eqref{ode} has a solution on $[0,1]$ and that $\delta(1)$ is integrable.
Notice that $q$ is twice continuously differentiable on $(l,r)$. The second
derivative
\[
q''(x)= \frac{2}{\eta^2(x)}
\]
is positive, which means that $q$ is convex. Moreover, $q$ is decreasing
on $(l,m)$ and increasing on $(m,r)$; in particular $q \ge0$.

%
\begin{theorem}\label{thm:intcondidelta}
If the function $q$ is integrable wrt $\nu$, then the ODE \eqref{ode}
has a
solution on $[0,1]$ for almost all paths of $M$ and $\delta(1)$ is integrable.
In this case, $\EE[\delta(1)] = \int q(x) \nu(\mathrm{d}x)$.
\end{theorem}

\begin{pf}
Assume first that $q$ is integrable wrt $\nu$. This means that the random
variable $q(N_1)$ is integrable. Let
%
\begin{equation}
\label{defi tau_n} \tau_n = 1 \wedge\inf\biggl\{t \ge0\Bigm|\int
_0^t \bigl |q'(N_s)
b_x\bigl(s,B(s,N_s)\bigr)\bigr |^2 \,\mathrm{d}s
\geq n\biggr\},
\end{equation}
and observe that $(N_u)_{u \leq s}$ is bounded away from $l$ and $r$
for any
$s<1$,
and hence
$\tau_n \uparrow1$, a.s.
By It\^o's formula, and using $q(N_0)=q(m) =0$, we get
\begin{eqnarray*}
q(N_{\tau_n}) &=& \int_0^{\tau_n}
q'(N_s) b_x\bigl(s,B(s,N_s)
\bigr) \,\dd \Wtil_s + \frac{1}{2} \int_0^{\tau_n}
q''(N_s) b_x
\bigl(s,B(s,N_s)\bigr)^2 \,\dd s
\\
&=&\int_0^{\tau_n} q'(N_s)
b_x\bigl(s,B(s,N_s)\bigr) \,\dd\Wtil_s + \int
_0^{\tau_n} \frac{b_x(s,B(s,N_s))^2}{\eta^2(N_s)} \,\dd s.
\end{eqnarray*}
Taking expectations, the martingale part disappears and we obtain
%
\begin{equation}
\label{eq aux 160612} \EE \biggl[ \int_0^{\tau_n}
\frac{b_x(s,B(s,N_s))^2}{\eta^2(N_s)} \,\dd s \biggr] = \EE\bigl[q(N_{\tau_n})\bigr].
\end{equation}
Notice that Jensen's inequality implies that
\[
0 \le q(N_{\tau_n}) \le\EE\bigl[q(N_1) \mid
\Ftil_{\tau_n}\bigr].
\]
Since the family $(\EE[q(N_1) \mid \Ftil_{\tau_n}])_{n \ge1}$ is uniformly
integrable, also $(q(N_{\tau_n}))_{n \ge1}$ is uniformly integrable.
Therefore we can interchange the expectation operator and the limit $n
\to
\infty$ on the RHS of \eqref{eq aux 160612}. By monotone convergence,
we can do
so also on the LHS and hence we get
\[
\EE \biggl[ \int_0^1 \frac{b_x(s,B(s,N_s))^2}{\eta^2(N_s)} \,\dd
s \biggr] = \EE\bigl[q(N_1)\bigr] < \infty.
\]
Lemma~\ref{aux lem 160612} implies that the ODE \eqref{ode} has a
solution on
$[0,1]$ for almost all paths of $M$ and that $\delta(1)$ is integrable.
\end{pf}

The reverse statement of Theorem~\ref{thm:intcondidelta} holds true as-well,
that is, if $\delta(1)$ is integrable, then $q$ is integrable wrt $\nu
$. Indeed,
we next show that the existence of an \textit{arbitrary} integrable
solution of
the SEP implies that $q$ is integrable wrt $\nu$.

\begin{proposition}\label{thm:intimplintofq}
Any stopping time $\tau$ that solves the SEP satisfies
%
\begin{equation}
\label{expectation sol sep} \EE[\tau] \ge\int q(x) \,\mathrm{d}\nu(x).
\end{equation}
\end{proposition}

\begin{pf}
Recall that $\overline{\tau}= \tau\wedge H_l \wedge H_r$.
Since $\ell$ is absorbing if $H_{\ell}<\infty$ and similarly if $H_r
< \infty$
then $r$ is absorbing,
we have that
$M_{\overline{\tau}} = M_\tau$ and $\overline{\tau}$ is also an
embedding of
$\nu$.

Let $\tau$ be an stopping time with $M_{\tau} \sim\nu$. Suppose
that $\tau$
is integrable; else the statement is trivial. Let
\[
\sigma_n = n \wedge\inf\biggl\{t \ge0\Bigm| \int
_0^t \bigl\llvert q'(M_s)
\bigr\rrvert ^2 \eta(M_s)^2 \,\mathrm{d}s \ge
n\biggr\}.
\]

Observe that $\sigma_n \uparrow H_l \wedge H_r$, a.s. Using It\^{o}'s formula,
we obtain
%
\begin{eqnarray}
\label{ito q m stopped} \EE\bigl[q(M_{\tau\wedge\sigma_n})\bigr] &=& q(M_0) + \EE
\biggl[\frac{1}{2} \int_0^{\tau\wedge\sigma_n}
q''(M_s) \,\dd\langle M, M
\rangle_s \biggr]
\nonumber
\\[-8pt]
\\[-8pt]
&=& \EE[\tau\wedge\sigma_n].
\nonumber
\end{eqnarray}
Then Fatou's lemma implies
\[
\EE\bigl[q(M_\tau)\bigr]=\EE\bigl[q(M_{\overline{\tau}})\bigr] \le
\liminf_n \EE\bigl[ q(M_{\tau
\wedge
\sigma_n})\bigr] \le \EE[
\overline{\tau}] \leq\EE[\tau],
\]
and hence \eqref{expectation sol sep}.
\end{pf}

\begin{remark}
Notice that if $M$ attains the boundary $l$ with positive
probability in finite time, then the function $q$ is finite at $l$. In this
case $\nu$ can have mass on $l$. If $M$ does not attain the boundary
$l$ in
finite time, then obviously a distribution $\nu$ with mass in $l$ can
not be
embedded with an integrable stopping time. Similar considerations apply
at the
right boundary $r$.
\end{remark}

Theorems~\ref{thm:intcondidelta} and~\ref{thm:intimplintofq} imply the
following
corollaries.

\begin{coro}
Suppose $\nu\in L^1$ and $m = \int x \nu(\mathrm{d}x)$. There
exists an
integrable solution $\tau$ of the SEP
if and only
if $q$ is integrable wrt $\nu$. In this case, $\tau$ satisfies
\eqref{expectation sol sep}.
\end{coro}

\begin{coro}
Suppose $\nu\in L^1$ and $m = \int x \nu(\mathrm{d}x)$.
Whenever there exists an integrable solution of the SEP, then
$\delta(1)$ is also an integrable solution.
\end{coro}

\subsection{The non-centred case}

Suppose we are given a local martingale diffusion $M$ started at
$M_0=m$ and a
measure $\nu\in L^1$ with $\nu^* \neq m$.

Recall the definition of $q$ in (\ref{eqn:qdef}). To emphasise the
role of the
initial point, write $q_m$ for this expression. More generally, for $n
\in
(l,r)$ define\vspace*{-1pt}
%
\begin{equation}
\label{eqn:qdefG} q_n(x) = \int_n^x
\,\mathrm{d}y \int_n^y \,\mathrm{d}z
\frac{2}{\eta(z)^2}.
\end{equation}

Then $q_m(z) = q_n(z) + q_m(n) + q_m'(n)(z-n)$. As $q=q_m$, in particular
\[
q(z) = q_{\nu^*}(z) + q\bigl(\nu^*\bigr) + q'\bigl(\nu^*
\bigr) \bigl(z - \nu^*\bigr)
\]
{\spaceskip=0.185em plus 0.05em minus 0.02em and $\int q(z) \nu(\mathrm{d}z) = \int q_{\nu^*}(z) \nu(\mathrm{d}z) + q(\nu^*)$.
Hence,
$\int q(z) \nu(\mathrm{d}z)$ is finite if and only if $\int q_{\nu^*}(z)\* \nu
(\mathrm{d}z)$} is
finite.

We state the following theorem in the case $m> \nu^*$ which necessitates
$r=\infty$, and then $\ell\in[-\infty, m)$. There is a
corresponding result
for $m<\nu^*$ in which the condition $\lim_{n \uparrow\infty}
q(n)/n <
\infty$ is replaced by
$\lim_{n \uparrow\infty} q(-n)/n < \infty$.
Note that the limit $\lim_n q(n) / n$ is well defined because $q$ is convex.

\begin{theorem}
\label{thm:integrableNC} Suppose $m>\nu^*$.

Suppose $\int q(z) \nu(\mathrm{d}z)<\infty$ and $\lim q(n)/n < \infty$. Then
$\delta^*$ is an integrable embedding of $\nu$.

Conversely, suppose there exists an integrable embedding $\tau$ of
$\nu$ in
$M$.
Then $\int q(z) \nu(\mathrm{d}z)<\infty$ and $\lim q(n)/n < \infty$.
\end{theorem}

\begin{pf}
Consider the first part of the theorem.
By the comments before the theorem, we may assume that $\int q_{\nu^*}(z)
\nu(\mathrm{d}z) < \infty$ and hence, for $M$ started at $\nu^*$,
$\EE[\delta_{\nu^*}(1)]<\infty$. Then, it is
sufficient to show that $\EE[H^M_{\nu^*}]<\infty$. But
\begin{eqnarray*}
\EE\bigl[H^M_{\nu^*}\bigr] &=& \lim_{n \uparrow\infty}
\EE \bigl[H^M_{\nu^*} \wedge H^M_n
\bigr] = \lim_{n \uparrow\infty} \EE\bigl[q(M_{H^M_{\nu^*} \wedge H^M_n})\bigr]
\\
& = & q\bigl(\nu^*\bigr) \lim_{n
\uparrow\infty} \frac{n-m}{n - \nu^*} + \lim
_{n \uparrow\infty} q(n) \frac{m - \nu^*}{n - \nu^*}
\\
& = & q\bigl(\nu^*\bigr) + \bigl(m - \nu^*\bigr) \lim_{n \uparrow\infty}
\frac{q(n)}{n},
\end{eqnarray*}
which is finite under the assumption that $\lim q(n)/n < \infty$.

For the converse, suppose that $\tau$ is an integrable embedding.
Without loss
of generality, we may assume that $\tau$ is minimal; if not we may
replace it
with a smaller embedding which is also integrable.
Then
\begin{eqnarray*}
\int q(x) \nu(\mathrm{d}x) &=& \EE\Bigl[\liminf_{n \to\infty}
q(M_{\tau\wedge H_{-n}^M \wedge H_{n}^M})\Bigr] \leq \liminf_{n \to\infty} \EE
\bigl[q(M_{\tau\wedge H^M_{-n} \wedge H^M_{n}})\bigr]
\\
& = & \lim_{n \to\infty} \EE\bigl[ \tau\wedge H^M_{-n}
\wedge H^M_{n}\bigr]
\\
& = & \EE[\tau] < \infty.
\end{eqnarray*}

It remains to show that $\EE[H^M_{\nu^*}] < \infty$. Recall that
this is
equivalent to the condition $\lim_{n \uparrow\infty} q(n)/n < \infty$.


Recall that by the Dubins--Schwarz theorem
(Rogers and Williams \cite{RogersWilliams}, page 64) we can write $M_t =
\hat{W}_{C_t}$
for a ${\mathbb G}=({\mathcal G}_t)_{t \geq0}$-Brownian motion $\hat{W}$
where ${\mathcal G}_s = {\mathcal F}_{C^{-1}_s}$. Let $\sigma= C_\tau
$. Then
$M_\tau= \hat{W}_{C_\tau} = \hat{W}_{\sigma}$ and $\sigma$ embeds
$\nu$ in
$\hat{W}$.

Since $\sigma$ is a
minimal embedding of $\nu$ in $\hat{W}$, by Theorem~5 of Cox and
Hobson \cite{CoxHobson05}
%
\begin{equation}
\label{eqn:ch1} \lim_n n \PP\bigl[\sigma>
H^{\hat{W}}_{-n}\bigr] = 0.
\end{equation}
Moreover, by arguments in the proof of Lemma~11 of Cox and
Hobson \cite{CoxHobson05}, for any stopping time $\tilde{\sigma}
\leq\sigma$
\[
\EE \bigl[ |\hat{W}_{\tilde{\sigma}}| \bigr] \leq\EE \bigl[ |\hat{W}_{\sigma}| \bigr] =
\int|z| \nu(\mathrm{d}z).
\]
Hence, $(\hat{W}_{t \wedge\sigma})_{t \geq0}$ is bounded in $L^1$,
and then
by Theorem~1
of Az\'{e}ma \textit{et al.} \cite{AzemaGundyYor}, $(\hat{W}_{t \wedge\sigma
})_{t \geq
0} $ is uniformly integrable if and only if
$\lim_n n \PP[\sigma> H^{\hat{W}}_{-n} \wedge H^{\hat{W}}_{n}] = 0$.
Since $\nu$ is not centred and $(\hat{W}_{t \wedge\sigma})_{t \geq
0} $ is not UI, it follows from (\ref{eqn:ch1}) that $\lim_n n \PP
[\sigma>
H^{\hat{W}}_{n}] >0$.
But $( \sigma> H^{\hat{W}}_{n}) \equiv(\tau> H^M_n)$ so
$\limsup_n n \PP[\tau> H^M_{n}] > 0$.

Then
\begin{eqnarray*}
\EE[\tau] & = & \lim_n \EE\bigl[q(M_{\tau\wedge H^M_n})\bigr]
\\
& \geq& \lim_n \EE\bigl[q(n) ; {\tau> H^M_n}
\bigr]
\\
& = & \lim_n \biggl( \frac{q(n)}{n} n \PP\bigl[\tau>
H^M_n\bigr] \biggr)
\\
& \geq& \lim\frac{q(n)}{n} \cdot\lim\sup n \PP\bigl[\tau>
H^M_n\bigr].
\end{eqnarray*}
Then, if $\EE[\tau]<\infty$ it follows that $\lim\frac{q(n)}{n} <
\infty$ and
$\EE[H^M_{\nu^*}] < \infty$.
\end{pf}

Finally, we consider the case where $\nu\notin L^1$.

%
\begin{lemma}
Suppose $\nu\notin L^1$. If $\tau$ is an embedding of $\nu$, then
$\tau$ is not integrable.
\end{lemma}

\begin{pf}
Observe that $q(x)\geq0$ and that if $\nu\notin L^1$ then since $q$
is convex we must have
$\int q(x) \nu(\mathrm{d}x) = \infty$. Then if $\tau$ is an embedding of $\nu$
\begin{eqnarray*}
\EE[\tau] &=& %
\lim_{n \uparrow\infty} \EE\bigl[\tau\wedge
H^{q(M)}_n \wedge H^{q(M)}_{-n} \bigr] =
\lim_{n \uparrow\infty} \EE\bigl[q(M_{\tau\wedge
H^{q(M)}_n\wedge H^{q(M)}_{-n}})\bigr]
\\
&\geq& \EE\Bigl[ \liminf_{n \uparrow\infty} q(M_{\tau\wedge
H^{q(M)}_n\wedge H^{q(M)}_{-n} })\Bigr] =
\EE\bigl[q(M_\tau)\bigr] = \infty.
\end{eqnarray*}
\upqed\end{pf}

\section{Bounded time embedding}
\label{sec:bte}

\subsection{The centred case}

In this section, we analyze the question under which conditions we can
guarantee the stopping time $\delta(1)$ to be bounded, that is,
$\delta(1)
\leq T \in\Rp$. Let us first state a necessary condition which places
a lower bound on how little mass must be embedded in each a
neighbourhood of a point~$x$.

\begin{theorem}
\label{thm:bte}
Suppose that $\eta$ is locally bounded and denote by $\eta^*$ its
upper semicontinuous envelope.
If a distribution with distribution function $F$ can be embedded before
time $T > 0$,
then for all $x \in\R$ with $0<F(x)<1$ it must hold that
%
\begin{equation}
\label{nece condi} \limsup_{\eps\dto0} - \eps^2 \ln
\bigl(F(x+\eps) - F(x - \eps)\bigr) \le \frac{\pi^2}{8} T \eta^*(x)^2.
\end{equation}
\end{theorem}

\begin{pf}
Fix $x$ and suppose $t'$ is such that $M_{t'}=x$.

For $\epsilon>0$ define $B_\epsilon(x)= \{ y \mid |y-x|<\epsilon\}$ and
$\bar{\eta}(x,\epsilon)=\max\{\eta^*(z)\mid z \in\bar B_\epsilon
(x)\} $.
Note that on $t \geq t'$ the process $\tilde{M}$ which solves the SDE
$\dd\tilde{M}_t = \tilde{\eta}(M_t)
\,\mathrm{d}W _t$ where
\[
\tilde{\eta}(m) = \bigl(1_{\{m \in B_\epsilon(x)\}} \eta(m) + 1_{\{m \notin B_\epsilon(x)\}} \bar{
\eta}(m,\epsilon) \bigr)
\]
subject to $\tilde{M}_{t'} = M_{t'}=x$,
coincides with $M$ up to the first leaving time of $B_\epsilon(x)$.
Moreover, there exists a Brownian motion $\tilde{W}$ such that on $t
\geq t'$,
$\tilde{M}_t = \tilde{W}_{\Gamma_t}$, where $\Gamma(t) = \int_{t'}^t \tilde{\eta}(M_s)^2 \,\dd s \leq
\bar{\eta}(x,\epsilon)^2 t$.
Then
\begin{eqnarray*}
\PP \Bigl[ \sup_{t' \leq t \leq T} \llvert M_t -
M_{t'} \rrvert < \epsilon \Bigr] &=& \PP \Bigl[ \sup
_{t'\leq t \leq T} | \tilde{M}_t -\tilde {M}_{t'}| <
\epsilon \Bigr]
\\
&=& \PP \Bigl[\sup_{t' \leq t \leq T}\llvert \tilde{W}_{\Gamma(t)} -
\tilde{W}_{\Gamma(t')}\rrvert < \epsilon \Bigr]
\\
&\geq& \PP \Bigl[ \sup_{0 \leq s \leq\bar{\eta}(x,\epsilon)^2 T} | {W}_s| < \epsilon
\Bigr].
\end{eqnarray*}
The probability for the absolute value of the Brownian motion ${W}$ to
stay within the ball
$B_\epsilon(0)$ up
to time $K T\ge0$ is given by (see Section~5, Chapter X in Feller
\cite{feller2})
\begin{eqnarray*}
\PP \Bigl[\sup_{s \in[0,K^2 T]} |W_s| < \epsilon \Bigr] &=&
\frac
{4}{\pi} \sum_{n = 0}^\infty
\frac{1}{2n+1} \mathrm{e}^{-\fraca{(2n+1)^2 \pi^2}{(8
\epsilon^2)} K T} (-1)^{n}
\\
&\ge& \frac{4}{\pi} \mathrm{e}^{-\fraca{\pi^2}{(8 \epsilon^2)} K T} - \frac
{4}{3 \pi}
\mathrm{e}^{-\fraca{9 \pi^2}{(8 \epsilon^2)}
K T} \geq\frac{8}{3\pi} \mathrm{e}^{-\fraca{\pi^2}{(8 \epsilon^2)} K T}.
\end{eqnarray*}

Assume that there
exists a stopping time $\tau$ such that $M_{\tau}$ has the
distribution $F$. Denote
by $\zeta= \inf\{t \ge0\dvtx  M_t = x\}$ the first time the process $M$
hits $x$. Since $F(x) \notin\{0, 1\}$,
the event $A = \{ \zeta< \tau\}$ occurs with positive probability.

Let $\FF_{\zeta}$ be the $\sigma$-field generated by $M$ up to
time
$\zeta$ and observe that $A \in\FF_{\zeta}$. Note further
that the process $Z=(Z_h)_{h \geq0}$ given by $Z_h = M_{h + \zeta} -
M_\zeta$ is independent of
$\FF_{\zeta}$.

Now suppose $\tau$ is bounded by $T$.
The mass of $F$ on the ball $B_\eps(x)$ has to be at least as large as
the probability that $A$ occurs
and that $X$ stays
within the ball $B_\eps(x)$ between $\zeta$ and $T$. Therefore,
\begin{eqnarray*}
F(x+\eps) - F(x - \eps) &\ge& \PP \Bigl[A \cap\Bigl\{\sup_{\zeta\le s
\le T}
|M_s - M_\zeta| < \eps\Bigr\} \Bigr]
\\
&=& \PP[A] \PP \Bigl[\sup_{\zeta\le s \le T} |M_s -
M_\zeta| < \eps \Bigr]
\\
&\geq& \PP[A] \frac{8}{3\pi} \mathrm{e}^{-\fraca{\pi^2}{(8 \epsilon^2)} \bar
{\eta}(x,\epsilon)^2 T}.
\end{eqnarray*}
Hence, we have
\[
-\eps^2 \ln \bigl(F(x+\eps) - F(x - \eps) \bigr) \le
\eps^2 \ln \frac{3\pi}{8 \PP[A]} + \frac{\pi^2}{8} \bar{\eta}(x,
\epsilon)^2 T,
\]
which implies the result.
\end{pf}

Now we turn to the converse, and sufficient conditions for these to
exist an embedding of $\nu$ in bounded
time. Suppose again that $\nu\in L^1$ and $M_0=m=\int x \nu(\mathrm{d}x)$.

Recall the definition of $r$ in \eqref{eq:defrR}.
The first result is an immediate corollary of Theorem~\ref{thm:odesol}.

\begin{coro}\label{suff condi bound}
If $\lambda(t,y)^2$ is bounded by $T \in\Rp$, for all $y \in\R$
and $t\in[0,1]$, then the stopping time
$\delta(1)$ is also bounded by $T$.
\end{coro}

\begin{proposition}
Assume that $F$ is absolutely continuous and has compact support.
Suppose $F$ has density $f$.
If $\eta$ and $f$ are bounded away from zero, then the stopping time
$\delta(1)$ is bounded.
\end{proposition}

\begin{pf}
Note that $h' = \frac{\phi}{f \circ F^{-1} \circ\Phi}$ and
thus it follows from $f$ bounded away from zero
that $h'$ is bounded. Hence, $b_x$ is bounded and thus $\lambda(t,y)$
is bounded.
\end{pf}

\begin{lemma}
\label{etaconcave}
Suppose that $\eta$ is concave on $(l,r)$.
Let $F$ be an absolutely continuous distribution with\vspace*{2pt} $\supp(F)
\subseteq[l,r]$ and suppose that
$\sup_{x \in[l,r]} \frac{h'(x)}{\eta(h(x))} \leq\sqrt{T} < \infty
$. Then $F$
is embeddable in bounded time, and there exists an embedding $\tau$
with $\tau\leq T$.
\end{lemma}

\begin{pf}
We have $b(t,x)= (\phi_{1-t} \star h)(x)$ and
\[
b_x(t,x)= \bigl(\phi_{1-t} \star h'\bigr) (x)
\leq\sqrt{T} \bigl(\phi_{1-t} \star(\eta\circ h)\bigr) (x) \leq\sqrt{T}
\eta\circ(\phi_{1-t} \star h) (x) = \sqrt{T} \eta\bigl(b(t,x)\bigr)
\]
and then $\lambda(t,x)^2 \leq T$ and the result follows from Corollary~\ref{suff condi bound}.
\end{pf}

\begin{remark}
\label{rem:notconcave}
More generally for the existence of a bounded embedding it is
sufficient that
there is a
concave function
$\xi$ and $\epsilon$ in $(0,1)$ for which $\epsilon\xi\leq\eta
\leq
\epsilon^{-1} \xi$. Then if $\sup_{x \in[l,r]} \frac{h'(x)}{\eta
(h(x))} \leq\sqrt{T}$
\[
b_x(t,x) \leq\sqrt{T} \epsilon^{-1} \bigl(
\phi_{1-t} \star(\xi\circ h)\bigr) (x) \leq \sqrt{T}\epsilon^{-1}
\xi\circ (\phi_{1-t} \star h) (x) \leq \sqrt{T} \epsilon^{-2}
\eta\bigl(b(t,x)\bigr),
\]
and $\lambda(t,x)^2 \leq T \epsilon^{-4}$.
\end{remark}

%
%
\begin{remark}
The sufficient condition from Lemma~\ref{etaconcave} implies a
stronger version of the necessary condition of Theorem~\ref{thm:bte}.
Indeed, it implies that the limit superior of the left hand side of
equation \eqref{nece condi} is equal to zero. To show this, let $\eta
$ be locally bounded and assume that $F$ satisfies the assumptions of
Lemma~\ref{etaconcave}. Then for all $z \in(\underline v, \bar v)$,
the interior of $\supp(F)$, we have
\[
f(z) \ge\frac{1}{\sqrt{T}} \frac{\phi\circ\Phi^{-1} \circ
F(z)}{\eta(z)}.
\]
Let $x \in(\underline v, \bar v)$. Since $\eta$ is locally bounded
there exists $B \in\R_+$ such that $\eta(z) \le B$ for $z$ close
enough to $x$. Then, for $\eps$ small we have
\begin{eqnarray*}
\ln\bigl(F(x+\eps) - F(x - \eps)\bigr) &=& \ln \biggl( \int_{x-\eps}^{x +
\eps}
f(z)\,\mathrm{d}z \biggr)
\\
&\ge& \ln \biggl( \frac{1}{\sqrt{T}B}\int_{x-\eps}^{x + \eps}
\phi\circ\Phi^{-1} \circ F(z)\,\mathrm{d}z \biggr),
\end{eqnarray*}
and applying Jensen's inequality we obtain
\[
\ln\bigl(F(x+\eps) - F(x - \eps)\bigr) \ge \ln\frac{2\eps}{\sqrt{T}B} -
\frac{1}{4\eps} \int_{x-\eps}^{x + \eps} \bigl[\bigl(
\Phi^{-1} \circ F(z)\bigr)^2 + \ln2 \pi\bigr] \,\mathrm{d}z.
\]
Consequently $\limsup_{\eps\downarrow0} - \eps^2 \ln(F(x+\eps) -
F(x - \eps)) = 0$.
\end{remark}

%
%
\subsection{The non-centred case}
\label{ss:Bnc}

If $M$ is a martingale and
if $\tau$ is bounded by $L$, then $M_{t \wedge\tau}$ is uniformly integrable
and $\EE[M_\tau]=m$. Hence, there are no embeddings of $\nu$ in $M$
if $\nu
\notin L^1$ or $\nu^* \neq m$.

If $M$ is a local martingale but not a martingale, then we may have
$\tau\leq
L$ and $\EE[M_\tau] \neq M_0 = m$. However, $\delta^*$ is not
bounded since
$\delta^* > H^M_{\nu^*}$ which is not bounded.

For example, let $m=1$ and $\eta(x)=x^2$ so that $\mathrm{d}M _t= M_t^2 \,\mathrm{d}B_t$,
and $M$
is the reciprocal of a 3-dimensional Bessel process. Let $\nu=
{\mathcal
L}(M_1)$. Then $\nu\in L^1$ and $\nu^* \leq1 = m$. Then, trivially,
$\tau
\equiv1 $ is a bounded embedding of $1$.

\section{General diffusions}
\label{sec:diff}

Let $(X_t)_{t \geq0}$ be a solution to
\[
\mathrm{d} X_t = \beta(X_t)\,\mathrm{d}t +
\alpha(X_t)\,\mathrm{d}W _t, \qquad\mbox{with }
X_0 = x_0,
\]
where $x_0 \in\R$, $\beta\dvtx  \R\rightarrow\R$ and $\alpha\dvtx  \R
\rightarrow\R$ are
Borel-measurable. We assume that $X$ takes values only in an interval
$[l,r]$ with
$-\infty\le l < x_0 < r \le\infty$. Moreover, we assume that $\alpha
(x) \neq 0$ for
all $x \in(l,r)$ and that $\frac{1+ |\beta|}{\alpha^2}$ is locally
integrable on
$(l,r)$.

Suppose we want to embed $\rho$ in $X$ with a stopping time $\tau$.

By changing the space scale one can transform the diffusion $X$ into a
continuous local martingale.
To this end, we define the \textit{scale function} $s$ (cf. \cite{RY},
Chapter VII, \S3)
via
\[
s ( x ) = \int_{x_0}^x \exp \biggl( - \int
_{x_0}^y \frac{2 \beta
( z )}{\alpha( z )^2} \,\dd z \biggr) \,\dd
y, \qquad x \in(l,r).
\]
Note that we are always free to choose the scale function such that
$M_0=s(x_0)=0$, and we have done so.

Then $s$ solves $\beta(x) s'(x) + \frac{1}2 \alpha^2(x) s''(x) = 0$.
Note that the scale function $s$ is strictly increasing and
continuously differentiable.
It\^o's formula implies that $M_t = s ( X_t )$ is a local martingale with
integral representation
\[
M_t = \int_0^t s' (
X_s ) \alpha( X_s ) \,\dd W_s.
\]
Thus $\dd M_t = \eta(M_t) \,\dd W_t$ where $\eta\equiv(s' \alpha)
\circ s^{-1}$.

Note that
\[
\int_{\chi-\epsilon}^{\chi+\epsilon} \frac{1}{((s' \alpha)\circ
s^{-1})^2(z)} \,\mathrm{d}z =
\int_{s^{-1}(\chi-\epsilon)}^{s^{-1}(\chi+\epsilon)} \frac
{1}{\alpha^2(z)s'(z)} \,\mathrm{d}z.
\]
Since $s$ is continuous, $\frac{1}{\eta}$ is locally square
integrable provided
\[
\frac{1}{\alpha(y) \sqrt{s'(y)}} = \frac{1}{\alpha(y)} \exp \biggl( -\int_{x_0}^y
\frac{2 \beta(z)}{ \alpha(z)^2 } \,\mathrm{d}z \biggr)^{-\fraca{1}{2}} = \frac{1}{\alpha(y)} \exp
\biggl( \int_{x_0}^y \frac{\beta(z)}{ \alpha(z)^2 } \,
\mathrm{d}z \biggr)
\]
is locally square integrable which follows from our assumptions on the pair
$(\alpha, \beta)$.

Let $F_\rho$ be the distribution function of $\rho$. If $\nu= \rho
\circ s^{-1}$ so that
$F(x) = F_\rho(s^{-1}(x))$,
then $X_\tau\sim\rho$ is equivalent to $M_\tau\sim\nu$.
Then $\nu$ has mean zero if and only if
%
\begin{equation}
\label{meancondition} \int_\R s(x) \rho(\mathrm{d}x) = 0.
\end{equation}

Clearly the requirement that $\tau$ is finite, integrable or bounded is
invariant under the change of scale. However, in the case of bounded
embeddings we can give a simple sufficient condition in terms of data relating
to the general diffusion $X$.

Define $g = F_\rho^{-1} \circ\Phi$ and $h = s \circ g = F^{-1}
\circ\Phi$.

\begin{theorem}
\label{thm:gendiff}
If $ x \mapsto- \frac{2 \beta(x)}{\alpha(x)} + \alpha'(x) $ is
non-increasing
and $\frac{g'}{\alpha\circ g}$ is bounded by $\sqrt{T}$,  then $\rho
$ can be
embedded in $X$ in bounded time. In particular, there exists an embedding
$\tau$ with $\tau\leq T$.
\end{theorem}

\begin{pf}
We prove in the first step that $\eta\equiv(s' \alpha) \circ s^{-1}$
is concave. We have
\[
\eta' = \bigl( \bigl(s' \alpha\bigr) \circ
s^{-1} \bigr)' = \frac{(s'' \alpha+ s' \alpha') \circ s^{-1}}{s' \circ s^{-1}} = \biggl( -
\frac{2\beta}{\alpha} + \alpha' \biggr) \circ s^{-1}
\]
where we have used the fact that $s$ solves $\alpha s'' = - 2\beta
s'/\alpha$.
As $s^{-1}$ is monotone increasing, under the first hypothesis of the
theorem we have that $ ( (s' \alpha)
\circ
s^{-1}  )'$ is
non-increasing and hence $\eta$ is concave.

We have $h = s \circ g$ and hence, again by hypothesis,
\[
\frac{h'}{\eta\circ h} = \frac{(s' \circ g) g'}{(s'\alpha) \circ g
} = \frac{g'}{\alpha
\circ g} \le\sqrt{T}.
\]
%
Lemma~\ref{etaconcave} implies that $\nu$ can be embedded in $M$
with a stopping time $\tau$ satisfying $\tau\leq T$, and the same stopping
time embeds $\rho$ in $X$.
\end{pf}

\section{Examples}
\label{sec:eg}

\subsection{Brownian motion with drift}

Let $X$ be a Brownian motion with drift, that is,
\[
X_t = x_0 + \gamma t + \theta W_t,
\]
where $\gamma\in\R$, $\theta> 0$ and $x_0 = 0$.
The scale function equals
\[
s(x) = %
\cases{ \displaystyle\frac{1}{\kappa} \bigl(1-\exp( - \kappa x ) \bigr)
&\quad$\mbox{for } \kappa\neq0$,\vspace*{2pt}
\cr
x &\quad$\mbox{for } \kappa= 0$, }
\]
with $\kappa=\frac{2 \gamma}{\theta^2}$.
If $\kappa>0$ then $s({\mathbb R}) = (-\infty, 1/\kappa)$, whereas if
$\kappa<0$ then $s({\mathbb R}) = (1/\kappa, \infty)$. Then, if
$M=s(X)$ we
have $\mathrm{d}M _t = \theta(1 - \kappa M_t) \,\mathrm{d}W _t$, and $M$ is a martingale.

Suppose the aim is to embed $\rho$.
Let $F_\rho$ be the distribution
function of $\rho$ and write $\nu=\rho\circ s^{-1}$.
Since $\rho$ is a measure on ${\mathbb\R}$, $\nu$ is a measure on
$(l,r)$ and any embedding $\tau$ is finite.
Note that\vspace*{-1pt}
\[
\nu^* = \int s(x) \rho(\mathrm{d}x) = \frac{1}{\kappa} \biggl( 1 - \int
_\R \mathrm{e}^{- \kappa x} \rho(\mathrm{d}x) \biggr).
\]
Then, by Proposition~\ref{prop:mneqbarnu} there is an embedding of
$\nu$ if and only if one of the following conditions is satisfied
\begin{enumerate}[3.]
\item$\nu^* \geq0$ and $\kappa> 0$.
\item$\nu^* \leq0$ and $\kappa< 0$.
\item$\kappa= 0$.
\end{enumerate}
Condition 1 and 2 simplify to
$%
0 \leq\nu^* \kappa= 1 - \int_\R \mathrm{e}^{- \kappa x} \rho(\mathrm{d}x)
$
and hence $\int_\R \mathrm{e}^{- \kappa x} \rho(\mathrm{d}x) \leq1$
is necessary for the existence of an embedding if $\kappa\neq0$.

\subsubsection{The centred case}

Suppose $\int \mathrm{e}^{- \kappa x} \rho(\mathrm{d}x) = 1$. Then $\nu$ has zero mean.

\begin{proposition}
For $\kappa\neq0$ ($\kappa=0$) there exists an integrable stopping time
embedding $\rho$ into $X$ if and only if $x$ ($ x^2$) is integrable with
respect to $\rho$. In this case, any minimal and integrable stopping time
$\tau$ satisfies\vspace*{-1pt}
\[
\EE[\tau] = %
\cases{ \displaystyle\frac{1}{\gamma} \int x \rho(\mathrm{d}x) &
\quad$\mbox{for } \kappa\neq0$,\vspace*{2pt}
\cr
\displaystyle\frac{1}{v^2} \int x^2\rho(
\mathrm{d}x) &\quad$\mbox{for } \kappa= 0$. } %
\]
\end{proposition}

\begin{pf}
Note that $q$ is given by\vspace*{-1pt}
\[
q(x) = %
\cases{ - \displaystyle\frac{2}{\kappa\theta^2} \biggl( \frac{1}{\kappa} \ln
(1 - \kappa x ) + x \biggr) &\quad$\mbox{for } \kappa\neq0$,\vspace*{2pt}
\cr
\displaystyle\frac{x^2}{\theta^2} &\quad$\mbox{for } \kappa= 0$. } %
\]
Moreover,\vspace*{-1pt}
\[
\int q(x) \nu(\mathrm{d}x) = \int q\bigl(s(x)\bigr) \rho(\mathrm{d}x) =
\cases{ \displaystyle\frac{1}{\gamma} \int x \rho(\mathrm{d}x) &\quad$\mbox{for
} \kappa\neq0$,\vspace*{2pt}
\cr
\displaystyle\frac{1}{v^2} \int x^2\rho(\mathrm{d}x) &
\quad$\mbox{for } \kappa= 0$. } %
\]
The result
follows now
from Theorem~\ref{thm:intcondidelta} and Proposition~\ref{thm:intimplintofq}.\vspace*{-1pt}
\end{pf}

Finally, we consider sufficient conditions for there to exist a bounded
embedding.
It turns out that the embedding stopping time
$\delta(1)$ is bounded if $h=F^{-1}\circ\Phi=s \circ F_\rho^{-1}
\circ\Phi$ is Lipschitz
continuous
with parameter $L$. We can thus recover the sufficient condition from
Section~3.2 in \cite{as11}.

%
\begin{proposition}
Suppose that $F_\rho^{-1}\circ\Phi$ is Lipschitz
continuous
with Lipschitz constant $L \in\Rp$. Then there exists an embedding
$\tau$ of $\rho$ in $X$ such that $\tau\le
\frac{L^2}{\theta^2}$.
\end{proposition}

\begin{pf}
For this example, $x \mapsto- 2 \beta/\alpha+ \alpha'$ is the constant
map. Hence, the result follows from Theorem~\ref{thm:gendiff}.
\end{pf}

\subsubsection{The non-centred case}

Suppose $\int \mathrm{e}^{-\kappa x} \rho(\mathrm{d}x) < 1$. It is clear that $\delta^*$
is finite almost surely, but the arguments of Section~\ref{ss:Bnc}
show that
there can
be no embedding of $\rho$ which is bounded. Further, if $\int
\mathrm{e}^{-\kappa
x} \rho(\mathrm{d}x) < 1$ then it follows that $\nu\in L^1$ and that
$\nu^* \in(0, 1/\kappa)$ (or $-1/\kappa,0)$.

Now consider integrable embeddings. By Theorem~\ref{thm:integrableNC}, there
exists an
integrable embedding if and only if $\EE[H^M_{\nu^*}] < \infty$
and $\int q(x) \nu(\mathrm{d}x)<\infty$. But
\( \EE[H^M_{\nu^*}] = \EE[H^X_{s^{-1}(\nu^*)}] \)
and, since $X$ is drifting Brownian motion, provided $\sgn(z)=
\sgn(\kappa) = \sgn(\gamma)$, $X$ hits $z$ in finite mean time.
Hence, $\EE[H^M_{\nu^*}]<\infty$. Further
\begin{eqnarray*}
\int q(x)\nu(\mathrm{d}x) & = & \int q\bigl(s(x)\bigr) \rho(\mathrm{d}x)
\\
& = & - \frac{2}{\kappa\theta^2} \int \biggl[ \biggl( \frac
{1}{\kappa} \ln \bigl(1 -
\kappa s(x) \bigr) + s(x) \biggr) \biggr] \rho(\mathrm{d}x)
\\
& = & - \frac{2}{\kappa\theta^2} \int \biggl[ -x + \frac{1}{\kappa} \bigl(1 -
\mathrm{e}^{-\kappa x}\bigr) \biggr] \rho(\mathrm{d}x)
\\
& = & \int\frac{x}{\gamma} \rho(\mathrm{d}x) - \frac{1}{\gamma} \int
\frac{1}{\kappa} \bigl(1 - \mathrm{e}^{-\kappa x}\bigr) \rho(\mathrm{d}x)
\\
& = & \frac{1}{\gamma} \biggl( \int x \rho(\mathrm{d}x) - \nu^* \biggr).
\end{eqnarray*}
Hence, there is an integrable embedding if $\int x \rho(\mathrm{d}x) < \infty$ and
$\delta^*$ is integrable.

\subsection{Bessel process}

Let $R$ be the radial part of 3-dimensional Brownian motion so that $R$ solves
$\mathrm{d}R_t = \mathrm{d}B_t + R_t^{-1}\,\mathrm{d}t $
and suppose that $R_0=1$. Then the scale function is given by $s(r)= 1
- r^{-1}$, and we can embed any
distribution $\rho$ on $\Rp$ in $R$ provided $\int r^{-1} \rho(\mathrm{d}r)
\le1$ (see Proposition~\ref{prop:mneqbarnu}).

\subsubsection{The centred case}

Suppose that $\int r^{-1} \rho(\mathrm{d}r) = 1$. Then $\nu$ has zero mean.

\begin{proposition}
There exists an integrable stopping time that embeds $\rho$ into $R$
if and only if $\int r^2 \rho(\mathrm{d}r) < \infty$. In this case, any
minimal and integrable stopping time $\tau$ satisfies
$\EE[\tau] = - \frac{1}3 + \int\frac{1}3 r^2 \rho(\mathrm{d}r)$.
\end{proposition}

\begin{pf}
Note that $\eta(x) = (s' \alpha) \circ s^{-1}(x) = (1-x)^2$. Moreover,
\[
q(x) = - \frac{2}3 x + \frac{1}3 \frac{1}{(1-x)^2} -
\frac{1}3.
\]
Notice that $\int q(x) (\rho\circ s^{-1})(\mathrm{d}x) = \int ( \frac{1}3
r^2 + \frac{2}3 r^{-1} - 1 ) \rho(\mathrm{d}r)$, and hence the result
follows from Theorem~\ref{thm:intcondidelta} and Proposition~\ref{thm:intimplintofq}.
\end{pf}

By Remark~\ref{rem:felkot},
we have that $M_t = 1
-R_t^{-1}$ is not a martingale (this is the Johnson--Helms example of a
strict local
martingale). Further,
the map $r \mapsto- 2 \beta(r)/\alpha(r) + \alpha'(r) = - 2/ r$ is
increasing.

However, suppose we want to embed a target law $\rho$ in $R$ in
bounded time, where the support of $\rho$
is bounded away from both $0$ and $\infty$ by $\hat{l}$ and $\hat
{r}$, respectively. Let $\bar{l} =
s(\hat{l})$ and $\bar{r}=s(\hat{r})$. Let $\hat{R}$ be the stopped
Bessel process $\hat{R}_t = R_{t \wedge
H_{\hat{l}} \wedge H_{\hat{r}}}$ and let $\bar{M} = s(\hat{R})$.
Then $\bar{M}$ is a martingale, which is
absorbed at both $\bar{l}$ and $\bar{r}$. Then a necessary condition
for there to exist an embedding of
$\nu$ in $\bar{M}$ in bounded time is that $\nu$ has support in
$[\bar{l},\bar{r}]$ and
$\int_{\bar{l}}^{\bar{r}}
x \nu(\mathrm{d}x) = 0$. Hence, a necessary condition for it to be possible to
embed $\rho$ in $R$ in bounded time
is
that $\int_{\hat{l}}^{\hat{r}} r^{-1} \rho(\mathrm{d}r) = 1$. By Remark~\ref
{rem:notconcave}, a sufficient
condition is
that $\int_{\hat{l}}^{\hat{r}} r^{-1} \rho(\mathrm{d}r) = 1$ and $\log
F_\rho^{-1} \circ\Phi$ is Lipschitz
continuous.

\subsubsection{The non-centred case}

Suppose $\int r^{-1} \rho(\mathrm{d}r) < 1$. It is clear that $\delta^*$
is finite almost surely, but the arguments of Section~\ref{ss:Bnc}
show that there can be no embedding of $\rho$ which is bounded.

Consider integrable embeddings. By Theorem~\ref{thm:integrableNC}, there
exists an
integrable embedding if and only if $\lim_{n \to\infty} \frac
{q(n)}{n} < \infty$ and $\int q(x) \nu(\mathrm{d}x)<\infty$.
For the first part, we have that
\[
\lim_{n \to\infty} \frac{q(n)}{n}= \lim_{n \to\infty}-
\frac{2}3 + \frac{1}3 \frac{1}{(1-n)^2 n} - \frac{1}{3 n} =
-\frac{2}{3} < \infty.
\]
Furthermore,
\begin{eqnarray*}
\int q(x)\nu(\mathrm{d}x) & = & \int q\bigl(s(x)\bigr) \rho(\mathrm{d}x) = \int
_{\Rp} \biggl( \frac{1}3 r^2 +
\frac{2}3 r^{-1} - 1 \biggr) \rho(\mathrm{d}x)
\\
&=& \int_{\Rp} \biggl( \frac{1}3 r^2 -
\frac{1}3 - \frac{2}3 s(x) \biggr) \rho(\mathrm{d}x)
\\
&=& \frac{1}{3} \biggl(\int_{\Rp} r^2
\rho(\mathrm{d}x) - 1 - 2 \nu^* \biggr).
\end{eqnarray*}
Hence, there exists an integrable stopping time if $r^2$ is integrable
with respect to $\rho$.

\subsection{Ornstein--Uhlenbeck process}

Let $X$ be an Ornstein--Uhlenbeck solving the SDE
\[
\mathrm{d}X_t = \xi X_t\,\mathrm{d}t + \sigma\,
\mathrm{d}W _t,
\]
where $\xi\in\R$, $\sigma> 0$ and $X_0=0$. The scale function is
given by $s(x) = \int_0^x \mathrm{e}^{(-\fraca{\xi}{\sigma^2} y^2)}\,\mathrm{d}y $.

\subsubsection*{The centred case}

Let $\rho$ be a distribution with $\int s(x) \rho(\mathrm{d}x) = 0$. Then $\nu
$ has zero mean.

We next give sufficient conditions for $\rho$ to be embeddable in
bounded time. We need to distinguish between a positive and negative
mean reversion speed $\xi$.

Suppose first that $\xi> 0$, and the process is mean repelling. Then
the scale function
is bounded. In this case,
$- \frac{2
\beta(x)}{\alpha(x)} + \alpha'(x) = -
\frac{2\xi}{\sigma} x$ is decreasing. Moreover,\vspace*{2pt} for $g=F^{-1}_\rho
\circ\Phi$,
$\frac{g'}{\alpha\circ g} =
\frac{1}{\sigma} g'$. Therefore, by Theorem~\ref{thm:gendiff}, if
$g$ is Lipschitz
continuous with Lipschitz constant $L$, then there exists an embedding
that is bounded
by $\frac{L^2}{\sigma^2}$.

Suppose next that $\xi< 0$. Then the derivative of the scale function
satisfies $s'(x)
\ge1$, $x \in\R$. Moreover, $\eta(x) = (s' \alpha) \circ s^{-1}(x)
\ge\sigma$ and the
intensity of the time change satisfies $r^2(t,x) \le\frac{1}{\sigma
^2} b^2_x(t,x)$.
Therefore, if $h = s \circ F^{-1}_\rho\circ\Phi$ is Lipschitz
continuous with Lipschitz constant $L$, then there exists an embedding
that is bounded
by $\frac{L^2}{\sigma^2}$. (Note that $h' = (s' \circ g) g' \geq g'$
so that the
requirement that $h$ is Lipschitz is stronger than the requirement that
$g$ is
Lipschitz.)

Finally, suppose that $\xi=0$. Then the scale function is the identity
function, and
the arguments from each of the last two paragraphs apply and yield the
same sufficient
condition.

\section*{Acknowledgements}

We thank an anonymous referee for many helpful comments.
Stefan Ankirchner and Philipp Strack were supported by
the German Research Foundation (DFG) through the \textit{Hausdorff
Center for Mathematics} and SFB TR 15.



\printhistory


\begin{thebibliography}{19}


\bibitem{as11}
%
\begin{barticle}[mr]
\bauthor{\bsnm{Ankirchner},~\bfnm{Stefan}\binits{S.}} \AND
\bauthor{\bsnm{Strack},~\bfnm{Philipp}\binits{P.}}
(\byear{2011}).
\btitle{Skorokhod embeddings in bounded time}.
\bjournal{Stoch. Dyn.}
\bvolume{11}
\bpages{215--226}.
\bid{doi={10.1142/S0219493711003255}, issn={0219-4937}, mr={2836522}}
\end{barticle}
%
\bptok{imsref}%
\endbibitem

\bibitem{AzemaGundyYor}
%
\begin{bincollection}[mr]
\bauthor{\bsnm{Az{\'e}ma},~\bfnm{J.}\binits{J.}},
\bauthor{\bsnm{Gundy},~\bfnm{R.~F.}\binits{R.F.}} \AND
\bauthor{\bsnm{Yor},~\bfnm{M.}\binits{M.}}
(\byear{1980}).
\btitle{Sur l'int\'egrabilit\'e uniforme des martingales continues}.
In \bbooktitle{Seminar on {P}robability {XIV} ({P}aris, 1978/1979)
({F}rench)}.
\bseries{Lecture Notes in Math.}
\bvolume{784}
\bpages{53--61}.
\baddress{Berlin}: \bpublisher{Springer}.
\bid{mr={0580108}}
\end{bincollection}
%
\bptok{imsref}%
\endbibitem

\bibitem{Bass1983}
%
\begin{bincollection}[mr]
\bauthor{\bsnm{Bass},~\bfnm{Richard~F.}\binits{R.F.}}
(\byear{1983}).
\btitle{Skorokhod imbedding via stochastic integrals}.
In \bbooktitle{Seminar on Probability {XVII}}.
\bseries{Lecture Notes in Math.}
\bvolume{986}
\bpages{221--224}.
\baddress{Berlin}: \bpublisher{Springer}.
\bid{doi={10.1007/BFb0068318}, mr={0770414}}
\end{bincollection}
%
\bptok{imsref}%
\endbibitem

\bibitem{CarrLee}
%
\begin{barticle}[mr]
\bauthor{\bsnm{Carr},~\bfnm{Peter}\binits{P.}} \AND
\bauthor{\bsnm{Lee},~\bfnm{Roger}\binits{R.}}
(\byear{2010}).
\btitle{Hedging variance options on continuous semimartingales}.
\bjournal{Finance Stoch.}
\bvolume{14}
\bpages{179--207}.
\bid{doi={10.1007/s00780-009-0110-3}, issn={0949-2984}, mr={2607762}}
\end{barticle}
%
\bptok{imsref}%
\endbibitem

\bibitem{CoxHobson2004}
%
\begin{barticle}[mr]
\bauthor{\bsnm{Cox},~\bfnm{A.~M.~G.}\binits{A.M.G.}} \AND
\bauthor{\bsnm{Hobson},~\bfnm{D.~G.}\binits{D.G.}}
(\byear{2004}).
\btitle{An optimal {S}korokhod embedding for diffusions}.
\bjournal{Stochastic Process. Appl.}
\bvolume{111}
\bpages{17--39}.
\bid{doi={10.1016/j.spa.2004.01.003}, issn={0304-4149}, mr={2049567}}
\end{barticle}
%
\bptok{imsref}%
\endbibitem

\bibitem{CoxHobson05}
%
\begin{barticle}[mr]
\bauthor{\bsnm{Cox},~\bfnm{A.~M.~G.}\binits{A.M.G.}} \AND
\bauthor{\bsnm{Hobson},~\bfnm{D.~G.}\binits{D.G.}}
(\byear{2006}).
\btitle{Skorokhod embeddings, minimality and non-centred target distributions}.
\bjournal{Probab. Theory Related Fields}
\bvolume{135}
\bpages{395--414}.
\bid{doi={10.1007/s00440-005-0467-y}, issn={0178-8051}, mr={2240692}}
\bptnote{check year}%
\end{barticle}
%
\bptok{imsref}%
\endbibitem

\bibitem{CoxWang2011}
%
\begin{barticle}[mr]
\bauthor{\bsnm{Cox},~\bfnm{Alexander~M.~G.}\binits{A.M.G.}} \AND
\bauthor{\bsnm{Wang},~\bfnm{Jiajie}\binits{J.}}
(\byear{2013}).
\btitle{Root's barrier: Construction, optimality and applications to
variance options}.
\bjournal{Ann. Appl. Probab.}
\bvolume{23}
\bpages{859--894}.
\bid{doi={10.1214/12-AAP857}, issn={1050-5164}, mr={3076672}}
\bptnote{check year}%
\end{barticle}
%
\bptok{imsref}%
\endbibitem

\bibitem{feller2}
%
\begin{bbook}[mr]
\bauthor{\bsnm{Feller},~\bfnm{William}\binits{W.}}
(\byear{1971}).
\btitle{An Introduction to Probability Theory and Its Applications.
{V}ol. {II}.},
\bedition{2nd} ed.
\blocation{New York}:
\bpublisher{Wiley}.
\bid{mr={0270403}}
\end{bbook}
%
\bptok{imsref}%
\endbibitem

\bibitem{Hobson2011}
%
\begin{bincollection}[mr]
\bauthor{\bsnm{Hobson},~\bfnm{David}\binits{D.}}
(\byear{2011}).
\btitle{The {S}korokhod embedding problem and model-independent bounds
for option prices}.
In \bbooktitle{Paris--{P}rinceton {L}ectures on {M}athematical
{F}inance 2010}.
\bseries{Lecture Notes in Math.}
\bvolume{2003}
\bpages{267--318}.
\baddress{Berlin}: \bpublisher{Springer}.
\bid{doi={10.1007/978-3-642-14660-2_4}, mr={2762363}}
\end{bincollection}
%
\bptok{imsref}%
\endbibitem

\bibitem{KS}
%
\begin{bbook}[mr]
\bauthor{\bsnm{Karatzas},~\bfnm{Ioannis}\binits{I.}} \AND
\bauthor{\bsnm{Shreve},~\bfnm{Steven~E.}\binits{S.E.}}
(\byear{1991}).
\btitle{Brownian Motion and Stochastic Calculus},
\bedition{2nd} ed.
\bseries{Graduate Texts in Mathematics}
\bvolume{113}.
\blocation{New York}:
\bpublisher{Springer}.
\bid{doi={10.1007/978-1-4612-0949-2}, mr={1121940}}
\end{bbook}
%
\bptok{imsref}%
\endbibitem

\bibitem{Kotani2006}
%
\begin{bincollection}[mr]
\bauthor{\bsnm{Kotani},~\bfnm{Shinichi}\binits{S.}}
(\byear{2006}).
\btitle{On a condition that one-dimensional diffusion processes are
martingales}.
In \bbooktitle{In Memoriam {P}aul-{A}ndr\'e {M}eyer: {S}\'eminaire de
{P}robabilit\'es {XXXIX}}.
\bseries{Lecture Notes in Math.}
\bvolume{1874}
\bpages{149--156}.
\baddress{Berlin}: \bpublisher{Springer}.
\bid{doi={10.1007/978-3-540-35513-7_12}, mr={2276894}}
\end{bincollection}
%
\bptok{imsref}%
\endbibitem

\bibitem{Monroe1972}
%
\begin{barticle}[mr]
\bauthor{\bsnm{Monroe},~\bfnm{Itrel}\binits{I.}}
(\byear{1972}).
\btitle{On embedding right continuous martingales in {B}rownian motion}.
\bjournal{Ann. Math. Statist.}
\bvolume{43}
\bpages{1293--1311}.
\bid{issn={0003-4851}, mr={0343354}}
\end{barticle}
%
\bptok{imsref}%
\endbibitem

\bibitem{OberhauserDosReis}
%
\begin{bmisc}[auto:STB|2014/02/12|14:17:21]
\bauthor{\bsnm{Oberhauser},~\bfnm{H.}\binits{H.}} \AND
\bauthor{\bsnm{dos Reis},~\bfnm{G.}\binits{G.}}
(\byear{2013}).
\bhowpublished{Root's barrier, viscosity solutions of obstacle
problems and reflected FBSDEs. Available at \arxivurl{arXiv:1301.3798}}.
\end{bmisc}
%
\bptok{imsref}%
\endbibitem

\bibitem{PedersenPeskir2001}
%
\begin{barticle}[mr]
\bauthor{\bsnm{Pedersen},~\bfnm{J.~L.}\binits{J.L.}} \AND
\bauthor{\bsnm{Peskir},~\bfnm{G.}\binits{G.}}
(\byear{2001}).
\btitle{The {A}z\'ema--{Y}or embedding in non-singular diffusions}.
\bjournal{Stochastic Process. Appl.}
\bvolume{96}
\bpages{305--312}.
\bid{doi={10.1016/S0304-4149(01)00120-X}, issn={0304-4149}, mr={1865760}}
\end{barticle}
%
\bptok{imsref}%
\endbibitem

\bibitem{RY}
%
\begin{bbook}[mr]
\bauthor{\bsnm{Revuz},~\bfnm{Daniel}\binits{D.}} \AND
\bauthor{\bsnm{Yor},~\bfnm{Marc}\binits{M.}}
(\byear{1999}).
\btitle{Continuous Martingales and {B}rownian Motion},
\bedition{3rd} ed.
\bseries{Grundlehren der Mathematischen Wissenschaften [Fundamental
Principles of Mathematical Sciences]}
\bvolume{293}.
\blocation{Berlin}:
\bpublisher{Springer}.
\bid{mr={1725357}}
\end{bbook}
%
\bptok{imsref}%
\endbibitem

\bibitem{RogersWilliams}
%
\begin{bbook}[mr]
\bauthor{\bsnm{Rogers},~\bfnm{L.~C.~G.}\binits{L.C.G.}} \AND
\bauthor{\bsnm{Williams},~\bfnm{David}\binits{D.}}
(\byear{1987}).
\btitle{Diffusions, {M}arkov Processes, and Martingales. {V}ol. 2}.
\bseries{Wiley Series in Probability and Mathematical Statistics:
Probability and Mathematical Statistics}.
\blocation{New York}:
\bpublisher{Wiley}.
\bid{mr={0921238}}
\end{bbook}
%
\bptok{imsref}%
\endbibitem

\bibitem{rogozin}
%
\begin{barticle}[mr]
\bauthor{\bsnm{Rogozin},~\bfnm{B.~A.}\binits{B.A.}}
(\byear{1966}).
\btitle{Distribution of certain functionals related to boundary value
problems for processes with independent increments}.
\bjournal{Teor. Veroyatn. Primen.}
\bvolume{11}
\bpages{656--670}.
\bid{issn={0040-361X}, mr={0208682}}
\end{barticle}
%
\bptok{imsref}%
\endbibitem

\bibitem{Root1969}
%
\begin{barticle}[mr]
\bauthor{\bsnm{Root},~\bfnm{D.~H.}\binits{D.H.}}
(\byear{1969}).
\btitle{The existence of certain stopping times on {B}rownian motion}.
\bjournal{Ann. Math. Statist.}
\bvolume{40}
\bpages{715--718}.
\bid{issn={0003-4851}, mr={0238394}}
\end{barticle}
%
\bptok{imsref}%
\endbibitem

\bibitem{SeelStrack2008}
%
\begin{barticle}[mr]
\bauthor{\bsnm{Seel},~\bfnm{Christian}\binits{C.}} \AND
\bauthor{\bsnm{Strack},~\bfnm{Philipp}\binits{P.}}
(\byear{2013}).
\btitle{Gambling in contests}.
\bjournal{J. Econom. Theory}
\bvolume{148}
\bpages{2033--2048}.
\bid{doi={10.1016/j.jet.2013.07.005}, issn={0022-0531}, mr={3146917}}
\end{barticle}
%
\bptok{imsref}%
\endbibitem


\end{thebibliography}
\end{document}